\documentclass[11 pt, a4paper]{article}

\usepackage{amsfonts,amsmath,amssymb,amsthm,url,float,enumitem}
\usepackage[pdftex]{graphicx}
\usepackage{pdfsync}
\usepackage{lscape}
\usepackage[hmargin=1.1in,vmargin=1.1in]{geometry}
\usepackage{color}
\usepackage{enumitem}
\usepackage{epsfig}
\usepackage{mathrsfs,dsfont}
\usepackage{caption}
\usepackage{hyperref}
\usepackage{jeffe_algo}

\newcommand{\ignore}[1]{}

\newtheorem{theorem}{Theorem}[section]

\newcommand{\parag}[1]{\vspace{2mm}

\noindent{\bf #1} }


\setlist{noitemsep}

\DeclareGraphicsExtensions{.jpg}
\definecolor{blue}{RGB}{65,105,225}

\newcommand{\pare}{\mathrm{par}}
\newcommand{\adj}{\mathrm{adj}}

\newcommand{\R}{\mathcal{R}}
\newcommand{\RG}{\mathcal{R}_G}

\title{$\R(K_6-e,K_4) =30$\thanks{This work was done as part of the 2022 Polymath Jr program, supported by NSF award DMS-2218374. This preprint has not undergone peer review or any post-submission improvements or corrections. The Version of Record of this article is published in Graphs and Combinatorics, and is available online at \url{https://doi.org/10.1007/s00373-025-02896-8}}}
\author{David James\thanks{University of Illinois Chicago, IL, USA. {\sl davidtj2@uic.edu}} 
\and Elisha Kahan\thanks{Yeshivah of Flatbush High School, NY, USA. {\sl elishakahan2.718@gmail.com}}
\and Erik Rauer\thanks{University of Minnesota-Morris, MN, USA. {\sl eriktr22@gmail.com}} }
\date{}

\begin{document}

\maketitle

\begin{abstract}
We settle the Ramsey problem $\R(K_6-e,K_4)$, also known as $\R(J_6,K_4)$ and $\R(K_6^-,K_4)$. 
Previously, the best bounds were $30\le \R(K_6-e,K_4) \le 32$.
We prove that $\R(K_6-e,K_4) =30$.
Our technique is based on the recent approach of Angeltveit and McKay and on older algorithms of McKay and Radziszowski.
\end{abstract}

\section{Introduction}

Ramsey theory suggests that every large object contains smaller structured pieces.
A classic example is that every red--blue edge coloring of the complete graph $K_6$ contains a red triangle or a blue triangle. 
For graphs $G_1,G_2$, let $\R(G_1,G_2)$ be the smallest integer $n$ such that every red--blue edge coloring of $K_n$ contains a red $G_1$ or a blue $G_2$.
The above example is part of the statement $\R(K_3,K_3)=6$.
Such expressions are called \emph{small Ramsey numbers}.

Discovering the exact value of small Ramsey numbers is quite challenging.
For example, while $\R(K_5,K_5)$ has attracted significant interest over many decades, we are far from knowing its exact value. 
The current best bounds are $43\le \R(K_5,K_5)\le 48$ (see \cite{AngelMcKay18,Exoo89}). 
An unusually large number of papers have been written about small Ramsey numbers.
A survey by Radziszowski on the subject \cite{Radziszowski11} is currently 116 pages long (without containing any proofs --- only problems and known results). 

Let $J_k$ be the graph on $k$ vertices with all possible edges except one.
In other words, $J_k$ is the complete graph $K_k$ with one edge removed.
This graph is also denoted as $K_k-e$ and as $K_k^-$.

In this work, we study the small Ramsey number $\R(J_6,K_4)$.
Recently, Lidicky and Pfender \cite{LP21} proved an upper bound of 32 for this number.
Boza \cite{Boza13} proved the lower bound 30. 
Thus, the best bounds were $30\le \R(J_6,K_4)\le 32$.
We settle the problem.

\begin{theorem} \label{th:main}
$\displaystyle \R(J_6,K_4) =30$.
\end{theorem}

Our basic approach follows the ideas of Angeltveit and McKay \cite{AngelMcKay18}.
We also rely on algorithms of McKay and Radziszowski \cite{MKS95}.
The proof is a mix of mathematical analysis and computations. 
Some of these computations use Python for simplicity and because the Python libraries NumPy and Dask provide good support for large arrays.
Other parts use Rust, to speed up the running time. 
For graph isomorphisms, we use nauty\footnote{See \url{http://users.cecs.anu.edu.au/~bdm/nauty/}}.

It seems plausible that a similar approach could lead to progress for similar problems, such as $\R(K_5,J_5)$, $\R(K_4,J_7)$, and $\R(J_4,J_8)$.
We may explore this direction in the future. 

Section \ref{sec:main} contains the main structure of the proof of Theorem \ref{th:main}.
Then, Sections \ref{sec:enumeration}--\ref{sec:SAT} contain the more technical and algorithmic aspects of the proof. 

Our data can be found in \url{https://geometrynyc.wixsite.com/ramsey}.
This includes the graphs that we enumerate and the code of all algorithms in this work. 

\parag{Notation.} Consider a graph $G$. Abusing notation, we also refer to the set of vertices of this graph as $G$. 
For example, the number of vertices in $G$ is $|G|$.
We may write $v\in G$ for a vertex $v$. 
Also, $G\setminus\{v\}$ refers to removing $v$ and the edges adjacent to it from $G$.

The \emph{complement} of a graph $G$, denoted $\overline{G}$, is a graph with the same vertex set as $G$.
An edge $e$ exists in $\overline{G}$ if and only if $e$ does not exist in $G$.
We note that $G$ contains an induced $\overline{K_i}$ if and only if $\overline{G}$ contains a $K_i$.

Instead of using colors to define a Ramsey problem, we use existing and non-existing edges. 
That is, $\R(G_1,G_2)$ is the minimal $n$ such that, for every graph $H$ with $n$ vertices, $H$ contains $G_1$ or $\overline{H}$ contains $\overline{G_2}$.
This is clearly equivalent to the red--blue approach. 
It simplifies some of our explanations below.  

Let $\RG(G_1,G_2)$ be the set of all graphs $H$ such that $H$ contains no $G_1$ and $\overline{H}$ contains no $\overline{G_2}$.
Let $\RG(G_1,G_2,m)$ be the set of all graph of $\RG(G_1,G_2)$ that have exactly $m$ vertices.

\section{Proof of Theorem \ref{th:main}} \label{sec:main}

This section consists of the general proof sketch of Theorem \ref{th:main}.
The more technical parts of the proof are deferred to later sections.

\begin{figure}[ht]
\centerline{\includegraphics[width=0.2\textwidth]{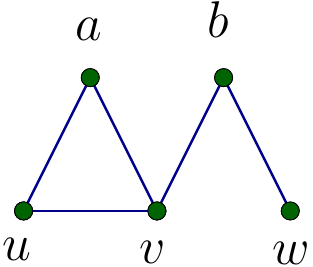}}
\caption{The subgraph $N(a)$ is induced by $u,v$. The subgraph $M(a)$ is induced by $b,w$. The subgraph $N(a,b)$ is the vertex $v$. The subgraph $N(a-b)$ is the vertex $u$.}
\label{fi:Neighbors}
\end{figure}

We assume, for contradiction, that there exists $F\in \RG(J_6,K_4,30)$.
For a vertex $a$ from $F$, let $N(a)$ be the subgraph of $F$ induced by the neighbors of $a$.
Let $M(a)$ be the subgraph induced by the vertices that are not neighbors of $a$.
Let $N(a,b)$ be the subgraph induced by the vertices that are neighbors of both $a$ and $b$.
Let $N(a-b)$ be the subgraph induced by the neighbors of $a$ that are not neighbors of $b$ and not $b$ itself.
See Figure \ref{fi:Neighbors}.
The vertex $a$ does not appear in $N(a),M(a),N(a,b)$, and $N(a-b)$.
The vertex $b$ does not appear in $N(a,b)$ and $N(a-b)$.

\parag{Vertex degrees.}
We claim that every vertex in $F$ has degree at least 13 and at most 18.
Indeed, consider a vertex $v$ from $F$. 
Since $\R(J_5, K_4) = 19$ (see \cite{EHM88}), if $\deg (v) \ge 19$ then $N(v)$ contains a $J_5$ or an induced $\overline{K_4}$.
The latter is a direct contradiction to $F\in \RG(J_6,K_4,30)$.
Considering $v$ together with a $J_5$ leads to a $J_6$, which is another contradiction.
Thus, all degrees in $F$ are at most 18. 
Since $R(J_6,K_3) = 17$ (see \cite{FRS80}), if $\deg (v)\le 12$ then $M(v)$ contains a $J_6$ or an induced $\overline{K_3}$.
A similar argument leads to a contradiction, so all degrees in $F$ are at least 13.

Since $F$ has 30 vertices and six possible degrees, there exists $13\le i\le 18$ such that $F$ contains at least five vertices of degree $i$.
It is impossible to have five vertices with each of the six degrees, since then the sum of the degrees in $F$ will be odd. 
Thus, there exists $13\le i\le 18$ such that at least six vertices of $F$ have degree $i$.
Since $F\in \RG(J_6,K_4,30)$, the six vertices of the same degree cannot form a $K_6$ or an induced $\overline{K_6}$.
We conclude that there is a pair of vertices of degree $i$ that are connected by an edge and another pair not connected by an edge.

\parag{The algorithm.}
For a fixed $i$ as defined above, let $a$ and $b$ be two vertices of degree $i$ that are connected by an edge.
That is, we have that $|N(a)|=|N(b)|=i$, that $b\in N(a)$, and that $a\in N(b)$.
We set $H=N(a,b)$ and $k = |H|$, which in turn implies that $|N(a-b)|=|N(b-a)|= i-k-1$.
We get that $N(a)\in \RG(J_5,K_4,i)$, since combining such a $J_5$ with $a$ leads to a $J_6$. 
Similarly, $N(b)\in \RG(J_5,K_4,i)$ and $N(a,b)\in \RG(J_4,K_4,k)$.
Since $H\in\RG(J_4,K_4)$ and $\R(J_4,K_4)=11$ (see \cite{ChH72}), we have that $k\le 10$.
Since $N(a-b)\in \RG(J_5,K_3)$ and $\R(J_5,K_3)=11$ (see \cite{Clancy77}), we have that $k\ge i-11$.
This implies that $N(a)\cap N(b)\neq \emptyset$.

\begin{figure}[ht]
\centerline{\includegraphics[width=0.43\textwidth]{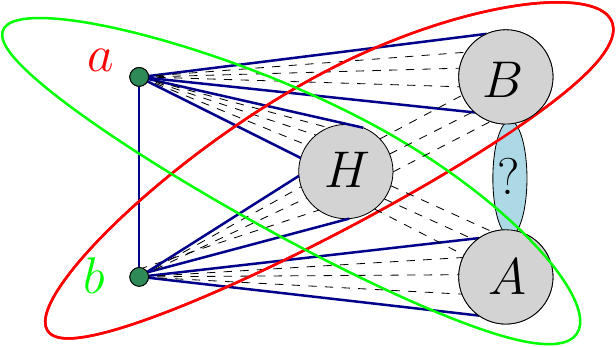}}
\caption{A big picture view of the analysis, following \cite{AngelMcKay18}. Here, $H=N(a,b)$, the green ellipse contains $N(b)$, the red ellipse contains $N(a)$, $A=N(b-a)$, and $B=N(a-b)$.}
\label{fi:BigPicture}
\end{figure}

A \emph{pointed graph} is a pair $(a,G)$ where $G$ is a graph and $a$ is a vertex of $G$. 
Our proof strategy is to enumerate all potential graphs $F\in \RG(J_6,K_4,30)$, as follows.
This is an inaccurate big picture strategy, for intuition.
Some details are changed later on.  

\begin{figure}[ht]
\centerline{\includegraphics[width=0.5\textwidth]{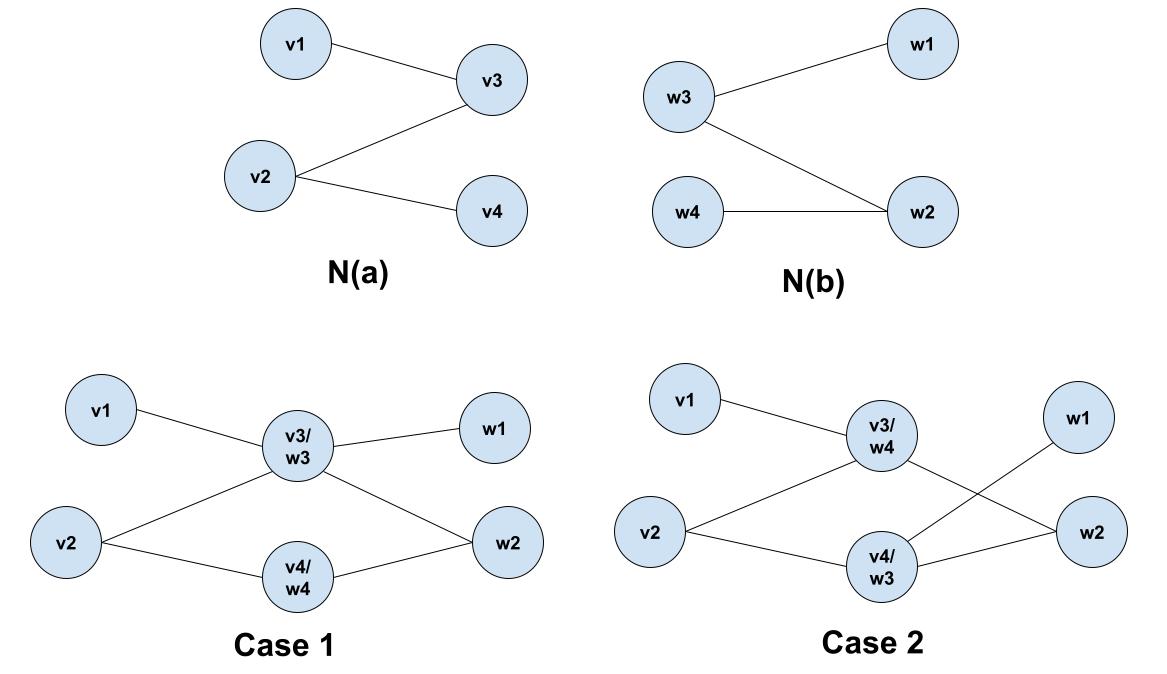}}
\caption{Each automorphism of $H$ leads to a different way of connecting $G_a$ and $G_b$.}
\label{fi:RamseyGluingMiniExample}
\end{figure}

In the following, $H=N(a,b)$, as in Figure \ref{fi:BigPicture}. 
Similarly, $G_a=N(a)$ and $G_b=N(b)$.
\begin{itemize}[noitemsep,topsep=1pt]
\item We enumerate all graphs of $\RG(J_4,K_4,k)$ and $\RG(J_5,K_4,i)$.
\item For each $H\in \RG(J_4,K_4,i)$ we consider all pairs of pointed graphs $(a,G_b),(b,G_a)$ such that $G_a,G_b\in\RG(J_5,K_4,i)$ and $H$ is an induced subgraph of both $G_b\setminus \{a\}$ and $G_a\setminus \{b\}$. 
We also ask for all vertices of $H$ to be connected to $a$ and $b$.
This leads to the graph in Figure \ref{fi:BigPicture}. 
Each automorphism of $H$ leads to a different way of connecting $G_a$ and $G_b$ at $H$. For example, see Figure \ref{fi:RamseyGluingMiniExample}.
\item By the above, $F$ contains at least one combination of $G_a,G_b,H$ as in the previous bullet. 
This combination may not be an \emph{induced} subgraph of $F$, since there might exist additional edges between the vertices of $A=N(a-b)$ and $B=N(b-a)$. 
Thus, we check all options for adding edges from $A\times B$, such that the resulting graph is in $\RG(J_6,K_4,2i-k)$. We refer to this process as \emph{gluing} $A$ and $B$.
\item For each graph $G$ generated above, we repeatedly check every way to add another vertex to $G$ while remaining in $\RG(J_6,K_4)$. We stop once no more vertices can be added. We refer to the process of adding a vertex as \emph{vertex extension}.
\end{itemize}

If the largest graph produced by the above algorithm contains fewer than 30 vertices, then we have a contradiction to the existence of $F$.
This contradiction implies that $\RG(J_6,K_4,30)$ is empty, so $\R(J_6,K_4)\le 30$.

Section \ref{sec:enumeration} describes the algorithm for enumerating the graphs of $\RG(J_5,K_4,i)$.
Section \ref{sec:glue} describes the gluing algorithm.
Section \ref{sec:VertexExt} describes the vertex extension algorithm.
Section \ref{sec:SAT} describes an alternative algorithm that performs the gluing and vertex extension simultaneously. 

\begin{figure}[ht]
\centerline{\includegraphics[width=0.6\textwidth]{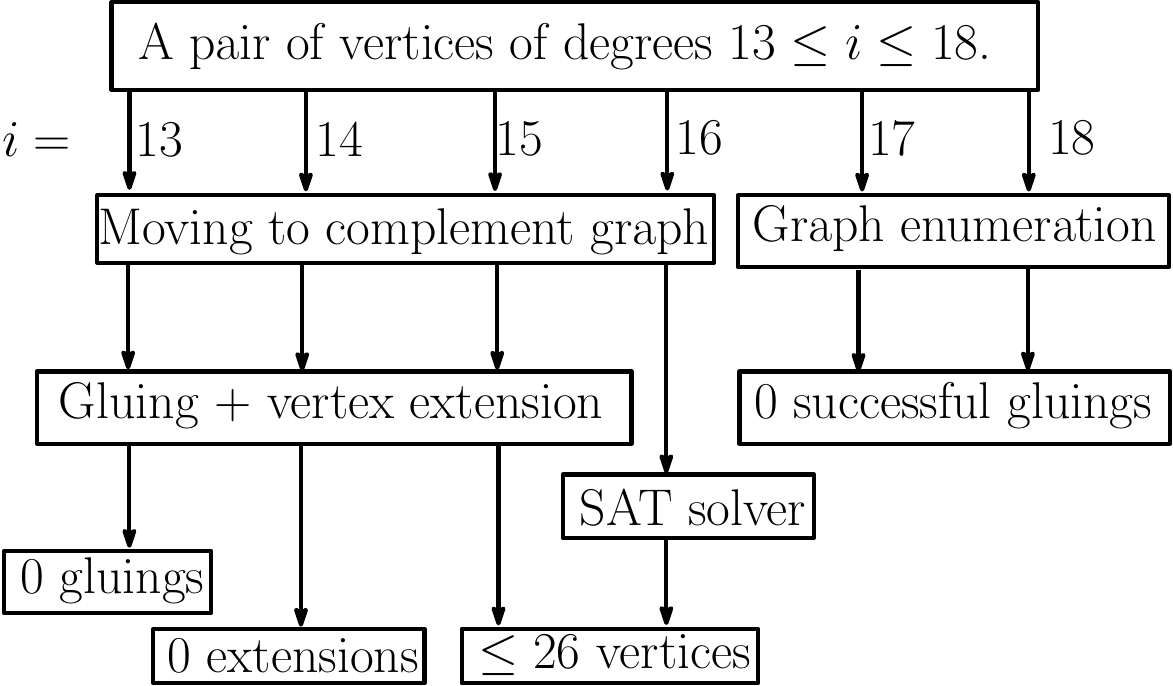}}
\caption{The flow of the case analysis.}
\label{fi:Flow}
\end{figure}

Recall that there exists $13\le i \le 18$ such that there is a pair of degree $i$ vertices that are connected by an edge and another pair not connected by an edge.
The rest of the analysis is divided into cases according to the value of $i$. 
While the above describes our general approach, some cases require modifications. 
Figure \ref{fi:Flow} depicts the flow of the case analysis.

\parag{The case analysis.} 
For the case of $i=18$, we enumerated the graphs of $\RG(J_5,K_4,18)$, obtaining six potential graphs for $N(a)$ and $N(b)$.
For $H$, the graphs of $\RG(J_4,K_4)$ are available at \cite{FidytekOnline1}.
Running the gluing algorithm does not lead to any successful gluings, so this case cannot occur.  
As a sanity check, we also ran the vertex extension algorithm on all six graphs of $\RG(J_5,K_4,18)$.
This leads to graphs with at most 24 vertices.

For the case of $i=17$, we enumerated the graphs of $\RG(J_5,K_4,17)$, obtaining $3,033$ potential graphs for $N(a)$ and $N(b)$.
For $H$, the graphs of $\RG(J_4,K_4)$ are available at \cite{FidytekOnline1}.
Running the gluing algorithm does not lead to any successful gluings, so this case cannot occur.  

Four cases remain: $13 \le i \le 16$.
In these cases, we take two vertices of degree $i$ with no edge between them and move to the complement graph $\overline{F}$. 
In $\overline{F}$, the two vertices are connected and of degree $29-i$. 
We note that $\overline{F}\in \RG(K_4,J_6)$, so $N(a),N(b)\in \RG(K_3,J_6)$ and $N(a,b)\in \RG(K_2,J_6)$.
Since $N(a,b)$ does not contain a $K_2$, it is an independent set. 
Since $\overline{N(a,b)}$ does not contain $J_6$, it has at most five vertices.

For the case of $i=16$, we have that $\deg(a)=\deg(b)=13$. 
The graphs of $\RG(K_3,J_6)$ are listed in \cite{FidytekOnline2}, and these are our $N(a)$ and $N(b)$.
Since $H$ is an independent set, it is easy to enumerate. 
The SAT solver algorithm from Section \ref{sec:SAT} leads to graphs with at most 26 vertices.
The computation took place over the course of approximately a month, using two home computers with the Intel Core i9-10900K CPU and Intel Core i7-13700K CPU, respectively. Not all cores were necessarily employed at once, and computations were manually shifted between the computers to share the workload.

For the case of $i=15$, we have that $\deg(a)=\deg(b)=14$. 
The graphs of $\RG(K_3,J_6)$ are listed in \cite{FidytekOnline2}, and these are our options for $N(a)$ and $N(b)$.
Since $H$ is an independent set, it is easy to enumerate. 
There are over seven billion successful gluings (before checking for isomorphisms). 
The vertex extension algorithm leads to graphs with 26 vertices, but not 27.

For the case of $i=14$, we have that $\deg(a)=\deg(b)=15$. 
The graphs of $\RG(K_3,J_6)$ are listed in \cite{FidytekOnline2}, and these are our options for $N(a)$ and $N(b)$.
Since $H$ is an independent set, it is easy to enumerate. 
After merging isomorphic gluing results, we obtain $1,477$ graphs of $\RG(K_4,J_6)$. 
The vertex extension algorithm fails for all these graphs.

For the case of $i=13$, we have that $\deg(a)=\deg(b)=16$. 
The graphs of $\RG(K_3,J_6)$ are listed in \cite{FidytekOnline2}, and these are our options for $N(a)$ and $N(b)$.
Since $H$ is an independent set, it is easy to enumerate. 
Running the gluing algorithm does not lead to any successful gluings, so this case cannot occur.  

Since all above cases lead to graphs with fewer than 30 vertices, we conclude that $\R(J_6,K_4)\le 30$.

\section{Graph Enumeration} \label{sec:enumeration}

In this section, we study the algorithm for enumerating the graphs of $\RG(J_5,K_4,\ell)$, for $\ell=17$ and $\ell=18$. 
Our general approach follows McKay and Radziszowski \cite{MKS95}, but with various changes.

For simplicity, we reverse the order, studying $\RG(K_4,J_5,\ell)$.
That is, we consider the complements of the graphs of $\RG(J_5,K_4,\ell)$.
Both $\RG(J_5,K_4,18)$ and $\RG(J_5,K_4,17)$ had been enumerated before, but we did not have access to these graphs. 
It is stated in the literature that $|\RG(J_5,K_4,18)|=6$ and $|\RG(J_5,K_4,17)|=3,033$ (for example, see \cite{Boza14}).
It seems that, a decade ago, some of these graphs were available online at \cite{Fidytek10}, but this is no longer the case.
We thus had to compute these graphs on our own.
Since we also received 6 and 3,033 graphs, the past results indicate that our enumeration is correct. 
We share our enumerated graphs at \url{https://geometrynyc.wixsite.com/ramsey}.

Consider $G\in \RG(K_4,J_5,\ell)$ and let $v$ be a vertex of $G$.
By definition, the subgraph induced by $N(v)$ is in $\RG(K_3,J_5)$ and the subgraph induced by $M(v)$ is in $\RG(K_4,J_4)$.
Thus, $G$ can be obtained by connecting a vertex $v$ to a $G_1\in \RG(K_3,J_5)$ and adding a $G_2\in \RG(K_4,J_4)$ that is not connected to $v$.
We set $m = |G_1|=\deg v$ and $m'=|G_2|=\ell -m-1$.
In Section \ref{sec:main} we proved that the degrees in a graph from $\R(J_6,K_4,30)$ are between 13 and 18.
Since $\R(K_3,J_5)=\R(K_4,J_4)=11$ (see \cite{ChH72,Clancy77}), the same analysis implies that $\ell-11 \le m \le 10$.
Since $m'=\ell-m-1$, we get that $\ell-11 \le m' \le 10$.

The graphs of $\RG(K_3,J_5)$ and $\RG(K_4,J_4)$ had been enumerated and are available online \cite{FidytekOnline1,FidytekOnline2}.
For the number of graphs of each type, see Table \ref{ta:enumeratedSizes}.

\begin{center}
\captionof{table}{The sizes of $\RG(K_3,J_5,i)$ and $\RG(K_4,J_4,i)$.}
\label{ta:enumeratedSizes}
\begin{tabular}{ ||c c c|| }
\hline
 $i$ & $|\RG(K_3,J_5,i)|$ & $|\RG(K_4,J_4,i)|$ \\ \hline \hline
 6 & 26 & 40 \\ \hline
 7 & 39 & 82 \\ \hline
 8 & 49 & 128 \\ \hline
 9 & 7 & 98 \\ \hline
 10 & 2 & 5 \\ \hline
\end{tabular} 
\end{center}

\begin{figure}[ht]
\centerline{\includegraphics[width=0.17\textwidth]{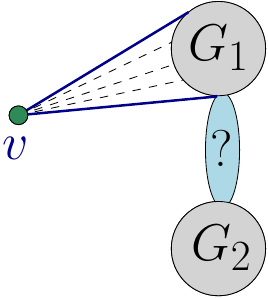}}
\caption{After fixing $G_1\in \RG(K_3,J_5)$ and $G_2\in \RG(K_4,J_4)$, it remains to choose the edges between $G_1$ and $G_2$.}
\label{fi:Enumeration}
\end{figure}

After combining $v,G_1$, and $G_2$, we need to decide which edges to add between $G_1$ and $G_2$. 
See Figure \ref{fi:Enumeration}.
Going over all possible edge choices and checking which are in $\RG(K_4,J_5,\ell)$ would take too long.
For example, when $\ell=18,m=8$, and $m'=9$, there are 72 potential edges in $G_1\times G_2$, so $2^{72}$ potential sets of edges.
By Table \ref{ta:enumeratedSizes}, in this case there are 49 options for $G_1$ and 98 options for $G_2$. 
Then, for each of the $2^{72}\cdot 49\cdot 98\approx 2\cdot 10^{25}$ resulting graphs, we need to check if it is in $\RG(K_4,J_5,\ell)$.

We use a more efficient approach to find the possible edge choices between $G_1$ and $G_2$.
This approach has two opposite directions: fixing a $G_1$ and finding all ways to glue it to all $G_2$ graphs, or fixing a $G_2$ and finding all ways to glue it to all $G_1$ graphs.
To optimize the running time, we usually fix a graph from the side with the fewer options.  
For example, when $m=9$, and $m'=8$, there are 7 options for $G_1$ and 128 options for $G_2$, so we fix $G_1$ and glue to it all 128 options for $G_2$.
We repeat this process for each of the 7 options for $G_1$.
When the graph numbers are similar, we fix the $G_2$ graph.
This case lends to better optimization with Dask due to its greater simplicity. 
The longest case involves gluing together graphs from $\R(K_4, J_4, 8)$ and $\R(K_3, J_5, 8)$, fixing graphs from the $\R(K_4, J_4)$ side, and taking about a week on a home computer using the Intel Core i9-10900K CPU.
We now describe both directions in detail.

\parag{Connecting a vertex of $G_1$ to vertices in $G_2$.}
In this case, we fix one graph $G_2$ and combine it with all possible options for $G_1$.
Denote the vertices of such a non-specific $G_1$ as $v_1,\ldots,v_{m}$.
A \emph{cone} of $v_i$ is a set of vertices of $G_2$ that we consider as a potential set of neighbors for $v_i$. 
A cone is \emph{feasible} if it does not lead to a $K_4$ in $G$ or to a $J_5$ in $\overline{G}$. 
Consider a feasible cone $C$ of $v_i$. 
Since $G_2$ does not contain $K_4$, the cone $C$ does not contain a $K_3$.
Since $\overline{G_2}$ does not contain a $J_4$ and $\R(K_3,J_4)=7$ (see \cite{ChH72}), every feasible cone consists of at most six vertices. 
Similarly, $G_2\setminus C$ does not contain an induced $\overline{K_3}$.  
Since $\R(K_4,K_3)=9$, the graph $G_2\setminus C$ has at most 8 vertices.

To enumerate all potential feasible cones, we go over all subgraphs of $G_2$ with at most six vertices.
For each such subgraph $C$, we keep it as a feasible cone if it contains no $K_3$ and $G_2\setminus C$ contains no induced $\overline{K_3}$ (by definition, $\overline{G_2}$ does not contain a $J_4$). 
This process is fast enough to be implemented in a straightforward way.
Denote the resulting feasible cones as $O_1,O_2,\ldots, O_n$.

Our next goal is to assign feasible cones to vertices of $G_1$.
We denote the cone assigned to $v_i$ as $C_i$.
We first define a set of rules for feasible cones, which do depend on the specific choice of $G_1$.
Below, we explain how the cone assignment is performed using these rules. 

\begin{itemize}[noitemsep,topsep=1pt] 
    \item[($k_2$)] Consider an edge $(v_i,v_j)$ from $G_1$. Then no edge has both of its points in $C_i\cap C_j$. Otherwise, we would have a $K_4$. 
    \item[($e_2$)] Consider vertices $v_i,v_j\in G_1$ that are \emph{not} connected by an edge. Then there is no $J_3$ in $\overline{G_2\setminus (C_i\cup C_j)}$. Otherwise, the complement will contain a $J_5$.
    \item[($e_3$)]Consider vertices $v_i,v_j,v_k\in G_1$ that form a $\overline{K_3}$. Then, for every edge $(u,u')$ in $G_2$, there exists at least one edge between $v_i,v_j,v_k$ and $u,u'$. For $u,u'$ that are not connected in $G_2$, there exist at least two edges between $v_i,v_j,v_k$ and $u,u'$.
    \item[($e_4$)] Consider an induced $\overline{K_4}$ in $G_1$. Then every $u$ in $G_2$ is in at least two cones of vertices of the $\overline{K_4}$.
    \item[($j_3$)] Consider a $J_3$ in $\overline{G_1}$. Then, for every $u,u'$ with no edge between them, at least one cone of a vertex from the $\overline{J_3}$ contains $u$ or $u'$.
    \item [($j_4$)] Consider a $J_4$ in $\overline{G_1}$. Then every vertex of $G_2$ is in at least one cone of a vertex from the $\overline{J_4}$.
\end{itemize}

It remains to explain how to use these rules to assign cones to vertices and how to simultaneously handle all graphs $G_1$.
We first discuss the algorithm of the other direction, and then describe the end of both algorithms together.

\parag{Connecting a vertex of $G_2$ to vertices in $G_1$.}
We start the analysis similarly to the start of the previous case. 
We fix one graph $G_1$ and combine it with all possible options for $G_2$.
We denote the vertices of such a non-specific $G_2$ as $u_1,\ldots,u_{m'}$.
A \emph{cone} of $u_i$ is a set of vertices from $G_1$ that are a potential set of neighbors for $u_i$. 
A cone $C'$ is \emph{feasible} when $\overline{G_1\setminus C'}$ does not contain a $J_4$ (by the definition of $G_1$, no cone contains a $K_3$ and no complement of a cone contains a $J_5$). 
Finally, $C'$ must contain at least two vertices of each induced $\overline{K_4}$ in $G_1$.
We enumerate all feasible cones $O'_1,\ldots, O'_n$, as before.

A feasible cone $C'$ of $u_i$ is \emph{minimal} if, when removing any vertex from $C'$, it is no longer feasible.
In other words, when removing any vertex from $C'$, we get that $\overline{G_1\setminus C'}$ contains a $J_4$, or that $C'$ contains a single vertex of an induced $\overline{K_4}$ in $G_1$.
We enumerate all minimal feasible cones by going over all cones, in increasing order of size. 
For each cone, we check if it is feasible and does not contain a minimal cone we already found. 
If these checks are successful, then we add the current set of vertices to our set of minimal cones. 

An \emph{interval} is a pair of feasible cones that are denoted \emph{top} and \emph{bottom}. 
We usually denote the top as $T$, the bottom as $B$, and the interval as $(B,T)$.
An interval $(B,T)$ must satisfy $B\subset T$.
We think of an interval $(B,T)$ as the set of all induced subgraphs of $H$ that contain all vertices of $B$ and do not contain any vertices not in $T$. 
To speed up our algorithm, we partition all feasible cones to disjoint intervals, as follows.

We create an ordered list $L$ of all feasible cones.
This list begins with the minimal cones in increasing order of size (number of vertices).
The non-minimal cones also appear in increasing order of size, after the minimal cones.
After creating $L$, we iterate through it.
When we reach a cone $B$ that is not part of an interval yet, we create a new interval with $B$ as its bottom.
The following paragraph explains how we find a top for this interval. 

To find a top $T$ for new bottom $B$, we first set $T$ to be the set of all vertices of $G_1$ (ignoring restrictions on the maximum cone size and being disjoint from other intervals).
We then go over each interval $(B',T')$ that was already created. 
If $B'$ contains a vertex not in $T$, then we move to check the next interval. 
If $B$ contains a vertex not in $T'$, then we move to check the next interval. 
Otherwise, for $(B,T)$ and $(B',T')$ to be disjoint, we choose a vertex of $B'$ and remove it from $T$. 
More specifically, the algorithm splits into different branches, each for removing a different vertex of $B'$ from $T$. 
Each such branch can split again when checking the following intervals. 
Once all branches are done, we take the largest $T$ to form a new interval with $B$.

The above process partitions all cones into disjoint intervals. 
Let $r$ be the number of intervals that we created.
As before, we require rules for interaction between different cones.
This time, instead of dealing with individual cones, the rules are about bottoms and tops of intervals.
Let the interval associated with $u_i$ be $(B_i,T_i)$.
\begin{itemize}[noitemsep,topsep=1pt]
    \item[($k'_2$)] Consider an edge $(u_i,u_j)$ from $G_2$. Then no edge in $G_1$ has both of its points in $B_i\cap B_j$. Also, $G_1\setminus (T_i \cup T_j)$ cannot contain an induced $\overline{K_3}$.
    \item[($k'_3$)] Consider a $K_3$ in $G_2$ with vertices $u_i,u_j,u_k$. Then $B_i \cap B_j \cap B_k$ is empty.  
    \item [($e'_2$)] Consider vertices $u_i,u_j$ from $G_2$ with no edge between them. Then $\overline{G_1\setminus(T_i\cup T_j)}$ cannot contain a $K_3$ or a $J_3$. Also, if $V$ is the set of vertices of an induced $\overline{K_3}$ in $G_1$ then $|T_i\cap V|+|T_j\cap V|>1$.
    \item[($e'_3)$] Consider an induced $\overline{K_3}$ in $G_2$ with vertices $u_i,u_j,u_k$. Then, $T_i\cup T_j \cup T_k=G_1$.     
    \item[($j'_3)$] Consider a $J_3$ in $\overline{G_2}$ with vertices $u_i,u_j,u_k$. Then there cannot be two vertices in $G_1\setminus (T_i\cup T_j \cup T_k)$ with no edge between them. 
\end{itemize} 

\parag{Combining rules with intervals.} 
We continue the process of assigning cones to vertices of $G_2$, by discussing how to apply the above rules to intervals, rather than to cones. 

We denote as $F(G_1, G_2, I_1, I_2,\ldots, I_{m'})$ the set of graphs where the feasible cones for $u_i$ are in $I_i$.
Let $I_i=(B_i,T_i)$.
If one of the above rules ($k'_2$), ($k'_3$), ($e'_2$), ($e'_3$), ($j'_3$) is violated, then there is no valid choice of cones for  $u_1,\ldots,u_m$.
However, when no rules are violated, there may still be bad cone choices. 
The following operations remove these bad choices.
\begin{itemize}[noitemsep,topsep=1pt]
    \item[($k'_2$)] Consider an edge $(u_i,u_j)$ from $G_2$. We remove from $T_i$ every $w\in B_j$ such that there exists an edge $(w,w')$ in $G_1$ with $w'\in B_i\cap B_j$. Also, we add to $B_i$ all vertices that are not in $T_j$ and form an induced $\overline{K_3}$ in $G_1$ with two vertices from $G_1\setminus (T_i \cup T_j)$.
    \item[($k'_3$)] Consider a $K_3$ in $G_2$ with vertices $v_i,v_j,v_k$. Then we remove from $T_i$ every vertex of $T_i \cap B_j \cap B_k$.  
    \item [($e'_2$)] Consider vertices $u_i,u_j$ from $G_2$ with no edge between them. We add to $B_i$ every vertex not in $T_j$ that froms a $J_3$ in the complement with two vertices from $G_1\setminus(T_i\cup T_j)$. 
We also add to $B_i$ every vertex that forms an induced $\overline{K_3}$ with two vertices from $G_1 \setminus (T_i \cup T_j)$.
    \item[($e'_3)$] Consider an induced $\overline{K_3}$ in $G_2$ with vertices $u_i,u_j,u_k$. We add to $B_i$ the vertices of $G_1\setminus (T_j \cup T_k)$.     
    \item[($j'_3)$] Consider a $J_3$ in $\overline{G_2}$ with vertices $u_i,u_j,u_k$. We add to $B_i$ the vertices of $G_1\setminus (T_j \cup T_k)$ that are not connected to another vertex from $G_1\setminus (T_i \cup T_j \cup T_k)$. 
\end{itemize} 

The above operations are not symmetric over $i,j,k$.
We thus apply each rule with each permutation of the relevant vertices.

Applying rules $(k'_2)$ and $(k'_3)$ may remove vertices from $T_i$ and thus lead to a new violation of the other three rules.
Similarly, applying rules $(e'_2),(e'_3),(j'_3)$ may add vertices to $B_i$ and thus lead to a new violation rules $(k'_2)$ and $(k'_e)$.
Let $F'(G_1, G_2, I_1, I_2,\ldots, I_{m'})$ be the result of repeatedly applying the above procedures until no interval needs to be revised. 
We say that $F'(G_1, G_2, I_1, I_2,\ldots, I_{m'})$ is \emph{collapsed}.
The \emph{collapsing process} is the process of repeatedly applying the above procedures until our objects are collapsed.
If each interval of $F'(G_1, G_2, I_1, I_2,\ldots, I_{m'})$ is a single cone, then it corresponds to a valid gluing. 
Otherwise, it may correspond to any number of valid gluings, including zero.

Being collapsed does not necessarily imply that all corresponding cone assignments are valid.  
For example, in rule $(k'_2)$ we remove from $T_i$ vertices that interfere with $B_j$, but we ignore vertices in cones larger than $B_j$.

\parag{Adding another tool for an improved running time.}
We now study the final part of the above algorithms for both cases. 
We start by explaining the second case, where we fix a $G_1$, since this case is more involved. 

\begin{figure}[ht]
\centerline{\includegraphics[width=0.6\textwidth]{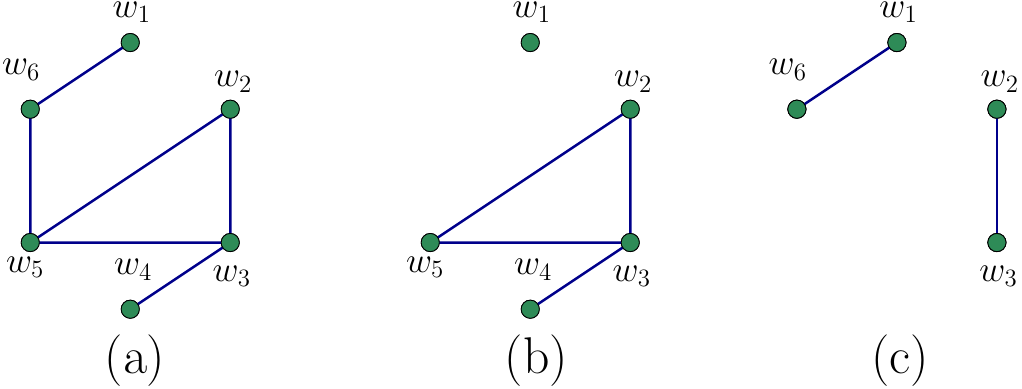}}
\caption{(a) A graph $G$. (b) The parent $\pare(G)$. (c) The adjunct $\adj(G)$ with respect to the sequence 1,1,2,3,4,4.}
\label{fi:ParAdj}
\end{figure}

Given a graph $G$ with vertices $w_1,\ldots,w_z$, the \emph{parent} of $G$ is the induced subgraph on $w_1,\ldots,w_{z-1}$. 
In other words, the parent is obtained by removing the last vertex with the edges adjacent to it.
See Figure \ref{fi:ParAdj}(a,b).
The \emph{adjunct} of $G$ is defined with respect to a sequence of integers $1= a_2 \leq a_3 \leq \cdots 
\le a_{z}$, where $a_i < i$.  
Two examples of valid sequences are $1, 2, 3, 4, 5$ and $1, 1, 1, 1, 1, 5, 5, 6$.
The adjunct of $G$ is the induced subgraph on $w_1,\ldots,w_{a_z - 1},w_z$. 
In other words, we remove the vertices $w_{a_z},w_{a_z+1},\ldots,w_{z-1}$.
See Figure \ref{fi:ParAdj}.
While the definition of an adjunct only relies on $a_{m'-1}$, we need the sequence to repeatedly perform the adjunct operation.

Let $\pare(G)$ and $\adj(G)$ denote the parent and adjunct of $G$, respectively.
We note that, if $F'(\pare(G_1),G_2, I_1, I_2,\ldots, I_{m'-1})$ or $F'(\adj(G_1),G_2, I_1, I_2,\ldots, I_{a_{m'-1}},I_{m'})$ have no valid gluings, then $F'(G_1, G_2, I_1, I_2,\ldots, I_{m'})$ cannot have valid gluings.

We now consider the case where both collapses $F'(\pare(G_1),G_2, I_1, I_2,\ldots, I_{m'-1})$ and $F'(\adj(G_1), G_2, I_1, I_2,\ldots, I_{a_{m'}-1},I_{m'})$ exist.
We denote the intervals that are produced by $F'(\pare(G_1), G_2, I_1, I_2,\ldots, I_{m'-1})$ as $I'_1,\ldots,I'_{m'-1}$.
We denote the intervals produced by $F'(\adj(G_1), G_2, I_1, I_2,\ldots, I_{a_{m'}-1},I_{m'})$ as $I''_1,\ldots,I''_{a_{m'}-1},I''_{m'}$. 
We will rely on the observation that $F'(G_1, G_2, I_1,\ldots,I_{m'})$ leads to the same gluings as
\[ F'(G_1, G_2, I'_1 \cap I''_1,\ldots,I'_{a_{m'} - 1} \cap I''_{a_{m'} - 1}, I'_{a_{m'}},\ldots, I'_{m' - 1}, I''_{m'}). \]
The above expression may not be fully collapsed, since it has not been checked if intervals from $I'_{a_{m'}},\ldots, I'_{m' - 1}$ violate any rules with $I''_{m'}$. 
Collapsing such violations may lead to additional rules being violated with other intervals. 

\begin{figure}[ht]
\centerline{\includegraphics[width=0.65\textwidth]{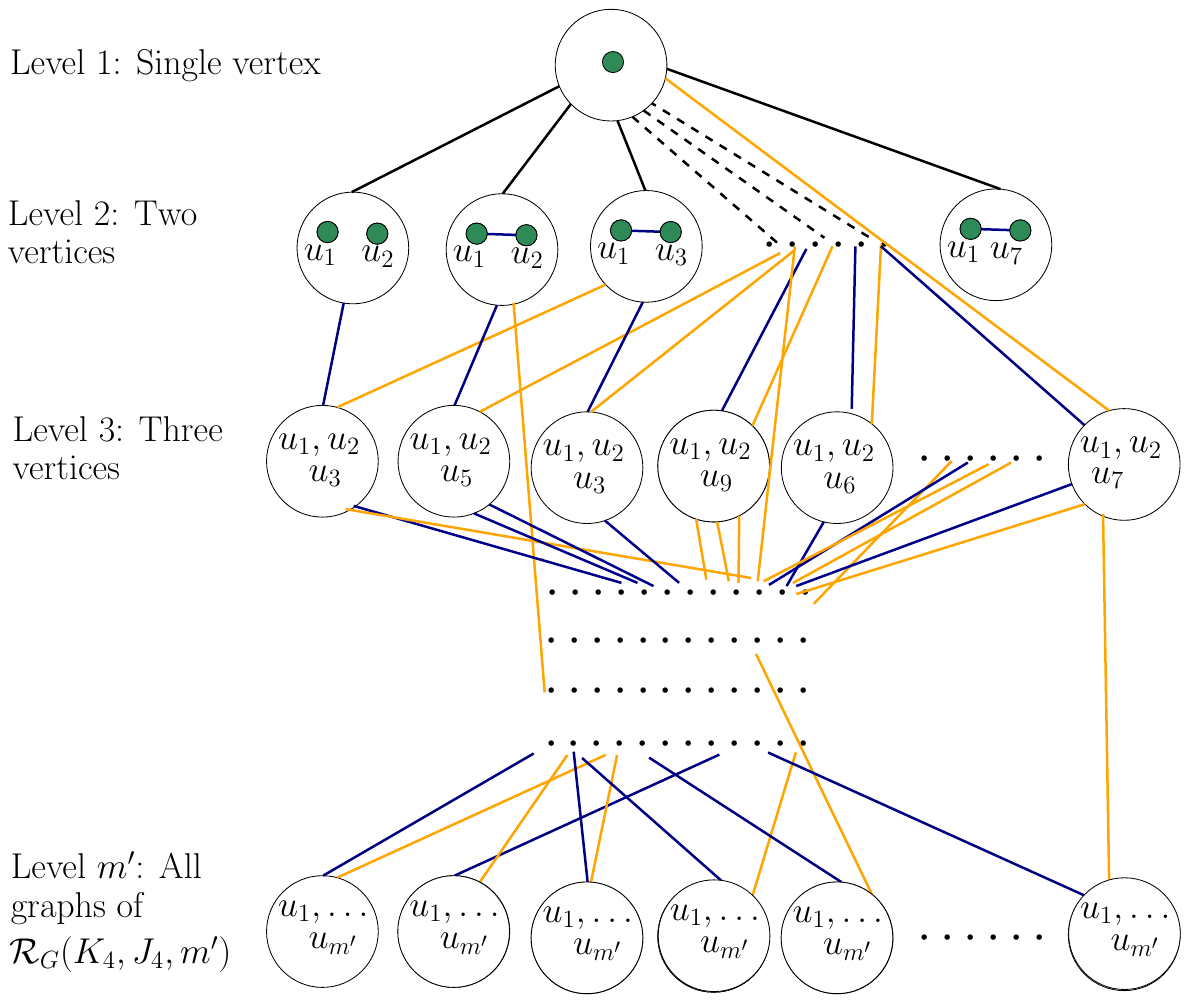}}
\caption{A double tree. Parent edges are blue and adjunct edges are orange.}
\label{fi:DoubleTree}
\end{figure}

A \emph{double tree} is a graph with two types of edges, which we denote as \emph{parent edges} and \emph{adjunct edges}.
When considering only the edges of any one type, the graph is a tree.
The flow of the improved algorithm is based on a double tree with $m'$ levels.
A node in level $i$ corresponds to a graph of $\RG(K_4,J_4,i)$.
Level $m'$ contains graphs of $\RG(K_4,J_4,m')$, which represent potential options for $G_2$. 
A tree node at level $i>1$ that corresponds to a graph $G$ is connected to two nodes in higher levels: a parent edge that connects the node to $\pare(G)$ and an adjunct edge that connects it to $\adj(G)$.
Note that $\pare(G)$ is at level $i-1$ and $\adj(G)$ can be at any level with an index smaller than $i$. 
See Figure \ref{fi:DoubleTree}.

It is not difficult to verify that the parent edges form a tree, and so do the adjunct edges.
A \emph{main branch} of the double tree starts at a level $m'$ node (a graph of $\mathcal{R}(K_4,J_4,m')$) and repeatedly travels up the tree, using only parent edges.
In Figure \ref{fi:DoubleTree}, a main branch is a path that uses only blue edges.

We are now ready to describe how the algorithm works.
We repeat the following for each $G_1 \in \RG(K_3,J_5)$.  
We create the feasible cones and intervals for this $G_1$, as described above. 
We then build the double tree, as follows.
The nodes of level $m'$ are the graphs of $\RG(K_4,J_4,m')$.
As stated above, these graphs are available at \cite{FidytekOnline2}, and we create a node for each.
We then iterate over every node, computing the parent and adjunct of the node and adding these new nodes and edges to the double tree.
When creating such a new node, we also add it to the set of nodes that are not processed yet. 
The \emph{root} of the double tree, at level 0, is a node corresponding to a graph with the single vertex.

The above process may generate multiple nodes with the same graph. 
For example, a node can be obtained in one way from a parent edge and in another way from an adjunct edge.
We check for isomorphisms and merge identical nodes, to keep the tree small.
Then, instead of a node containing an interval for each vertex from its graph, a node will contain an array where each cell holds an interval for each vertex.  
Since the value of $m'$ is not fixed, we build a separate double tree for each value of $m'$.
On the other hand, the same double tree can be used for all graphs $G_1$, so it suffices to build each double tree once.

Recall that adjuncts require a sequence $1= a_2 \le a_3 \le \cdots \le a_{z-1}$ with $a_i < i$.
We chose such sequences via experimentation --- checking which sequences make the algorithm run faster.
Intuitively, the larger $\adj(G)$ is, the longer it takes to collapse it. 
On the other hand, when $\adj(G)$ is larger, we expect smaller intersections between the intervals of the parent and the adjunct.
By the definition of parent edges and adjunct edges, every node in level $j>1$ corresponds to a graph with $v_1,\ldots,v_{j-1}$ and one additional vertex. 

When building the double tree, we also mark the nodes that belong to a main branch. 
This is easy to do: Whenever we process a node marked as being on a main branch, we also mark its parent as being on a main branch. 

Our end goal is to collapse every node on level $m'$, since this is equivalent to enumerating the graphs of $\RG(K_4,J_5,\ell)$.
We start at the root of the double tree and gradually travel down, handling the nodes of level $i$ before getting to level $i+1$.
However, we only collapse the level $i$ nodes that belong to a main branch. 
Recall that collapsing requires collapsed parent and adjunct. 
The parent is already collapsed by definition, but the adjunct might not be collapsed yet. 
If that is the case, we first collapse the adjunct, which might lead to more recursive collapsing. 

Recall that being collapsed does not imply that all corresponding cone assignments are valid.
In other words, the nodes of level $m'$ with no empty intervals include all graphs of $\RG(K_4,J_5,\ell)$, but possibly also other graphs. 
We thus continue the tree beyond level $m'$, as follows.
Consider a leaf node with at least one interval $(B,T)$ satisfying $B\neq T$.
For an arbitrary $w\in T\setminus B$, we create new child nodes where $(B,T)$ is respectively replaced with $(B\cup \{w\},T)$ and $(B,T\setminus\{w\})$.
We then collapse the two new child nodes and repeat the process for each.
This ends when each leaf of the tree contains an empty interval or corresponds to a single cone assignment.
We then check which of the latter type of leaves correspond to a graph of $\RG(K_4,J_5,\ell)$.

The above explains the double tree algorithm for the case where we fix a $G_1$ and simultaneously glue to it all options for $G_2$. 
We use the same algorithm when fixing a $G_2$ and simultaneously gluing all potential $G_1$ graphs to it. 
In that case, the algorithm is simpler, since there are no intervals. 
As before, each double tree node contains an array.
However, instead of intervals, each cell contains one cone for each vertex. 
In the collapsing process, we remove a cell if its cones violate one of the rules of this case.

\section{The Gluing Algorithm} \label{sec:glue}

In this section, we describe the algorithm for gluing the graphs $N(a)$ and $N(b)$, as mentioned in Section \ref{sec:main}. 
We only describe this algorithm briefly, since it is a variant of an algorithm of Angeltveit and McKay \cite[Section 5, second method]{AngelMcKay18}.

We work with the graph described in Figure \ref{fi:BigPicture}. 
Recall that \emph{gluing} is the process choosing edges between $A$ and $B$ without creating a $J_6$ or an induced $\overline{K_4}$. 
(If we are in the complement graph, we instead avoid $K_4$ and $J_6$ in the complement.)
We set $m=|A|=|B|=i-k-1$.
We denote the vertices of $A$ as $a_1,\ldots, a_m$, the vertices of $B$ as $b_1,\ldots,b_m$, and the vertices of $H$ as $h_1,\ldots,h_k$.

We create an $m\times m$ matrix $M$, where each cell contains one of the values True, False, or Unknown.
The value of cell $j$ in row $i$ states whether there is an edge between $a_i$ and $b_j$. 
At first, all matrix cells contain the value Unknown. 
Our goal is to change these values to True or False without creating copies of $J_6$ or induced $\overline{K_4}$.
A \emph{potential} $(r,s,t)$ \emph{set} is a set of $r$ vertices from $H$, $s$ vertices from $A$, and $t$ vertices from $B$. 
Such a potential set is a $j$\emph{-set} if $r+s+t=j$.
We only consider potential sets with $s>0$ and $t>0$, since sets with no pairs from $A\times B$ are unrelated to the gluing.

Consider a 6-set with $r+s=5$.
Such a set cannot contain a $J_6$, since that would imply that $N(b)$ contains a $J_5$. 
We may thus assume that $r+s\le 4$ and symmetrically that $r+t\le 4$. 
Since $r\ge 3$ implies $r+s\ge 5$ or $r+t\ge 5$, we conclude that it suffices to consider potential 6-sets with $r\le 2$.

We enumerate all potential 4-sets that have no edges between pairs of vertices not from $A\times B$ (no edges between two vertices from $H$, between a vertex from $H$ and a vertex from $A$, and so on). 
We will rely on these 4-sets to generate a gluing with no induced $\overline{K_4}$.
We also enumerate all potential 6-sets with at most one missing edge among pairs of vertices not from $A\times B$.
We will rely on these 6-sets to generate a gluing with no $J_6$.

We describe the algorithm in the original graph --- the complement case is symmetric.
We consider all potential $(2,1,1)$ sets: If the two vertices in $H$ are not connected to each other and to the other two vertices, then we set the cell of the edge between $A$ and $B$ to True.  
Otherwise, this 4-set will be an induced $\overline{K_4}$.
In the complement, the above argument for ignoring $J_6$ fails, so we check the 6-sets $(4,1,1)$, $(3,2,1)$, and $(3,1,2)$.

We create a stack $S$ and add to it all matrix cells that were set to True. 
We pop the top element $\alpha$ from $S$ and check each potential $(r,s,t)$ set that includes $\alpha$, as follows.
\begin{algo}
\textul{\textsc{ProcessStack}($M$ the $m \times m $ Matrix, $S$ the stack of matrix cells):}
\\While $S$ is not empty:\+
\\  Pop $\alpha$ the top element of $S$
\\  For each $(r,s,t)$ potential set $W$ that includes $\alpha$\+
\\      If $\alpha$ is False and $W$ is a 4-set with one Unknown edge $\beta$ and the other edges are False\+
\\          Change $M[\beta]$ to True and push $\beta$ to the stack $S$.\-
\\
\\      If $\alpha$ is True and $W$ is a 6-set with one edge False, one Unknown edge $\beta$, and the rest True\+
\\          Change $M[\beta]$ to True and push $\beta$ to the stack $S$.\-
\\
\\      If $\alpha$ is True and $W$ is a 6-set with two Unknowns, $\beta$ and $\gamma$, and the rest True\+
\\          Change both $M[\beta]$ and $M[\gamma]$ to False and push both to the stack $S$.\-
\\
\\      If the potential set leads to a forbidden configuration\+
\\          Declare that there are no valid gluings and stop.\-\-\-
\\ Check each $(r,s,t)$ potential set for forbidden configurations. 
\\Announce no valid gluings if any are found.
\end{algo}

If the above process ended by reaching to an empty $S$, we are not necessarily done, since there might still be Unknown edges.
In such a case, we arbitrarily choose an Unknown edge and split the process into two: one case where the edge is True and one where it is False.
We run the above process recursively for both cases. 
If we reach an empty $S$ and no Unknown edges, then this is a valid gluing to report. 

The performance of the described algorithm largely depends on the number of recursive calls. Usually, not many Unknowns are left before starting the recursive calls, so the above process runs in a reasonable time. After the recursive process ends, we add $a$ and $b$ to each resulting graph.

We briefly describe the differences of the algorithm in the complement case. Here, our above argument
for ignoring potential $6$-sets with $r \geq 3$ fails, so we check the 6-sets (4, 1, 1), (3, 2, 1), and (3, 1, 2). In the complement, we consider these sets at the start of the algorithm in addition to the $(2,1,1)$ sets in the expected ways to prevent $\overline{J_6}$'s from forming. Lastly, when we add $a$ and $b$ to our gluing at the end of our algorithm, this sometimes leads to copies of $K_4$ and of $J_6$ in the complement, so not all gluing results are valid.

\section{Vertex Extension} \label{sec:VertexExt}

In this section, we describe an algorithm that receives a graph $F\in \RG(K_4,J_6)$ and finds all ways of adding another vertex $w$ without creating a $K_4$ or a $J_6$ in the complement.
This is a variant of an algorithm from \cite{MKS95}.
We only need to find the set of neighbors of $w$. 
As in Section \ref{sec:enumeration}, we define an interval $I=[B,T]$ to represent all sets of vertices that contain $B$ and are contained in $T$.
In the current section, an interval represents possible sets of neighbors for $w$. 

We first enumerate all induced $K_3,\overline{J_5}$, and $\overline{K_5}$ in $F$.\footnote{This can be done using the Bron--Kerbosch algorithm \cite{BK73}.}
Let $X$ be the set of all such induced subgraphs.
We represent graphs as adjacency matrices.
Induced subgraphs are binary sequences, with a bit for each vertex. 
All steps of the following algorithm use bitwise operations, which lead to a fast running time.

We maintain a list $S$ of intervals that contain the neighbor sets we still consider. 
At first, since we have not disqualified any sets yet, $S$ contains one interval: $[\emptyset,F]$.
We then iterate over each element of $X$ and revise $S$ accordingly:


\begin{algo}
\textul{\textsc{VertexExtend(X)}:}
\\$S = [[\emptyset,F]]$
\\For $X_i \in X$:\+
\\  If $X_i$ is a $K_3$, For $I=[B,T] \in S$ with $X_i\subset T$:\+
\\    Remove $I$ from $S$; Label elements of $X_i\setminus B$ according to the sequence $v, v', v''$.
\\    If $|X_i\setminus B| \geq 1$, Add $[B,T\setminus \{v\}]$ to $S$. 
\\    If $|X_i\setminus B| \geq2$, Add $[B\cup\{v\},T\setminus \{v'\}]$.
\\    If $|X_i\setminus B| \geq 3$, Add $[B\cup \{v,v'\},T\setminus \{v''\}]$.\-
\\  If $X_i$ is a $\overline{J_5}$, For $I=[B,T] \in S$ with $X_i\cap B=\emptyset$:\+
\\    Remove $I$ from $S$; Label elements of $X_i\cap T$ according to the sequence $v, v', v'', \dots$
\\    If $|X_i\cap T| \geq 1$, Add $[B\cup \{v\},T]$ to $S$. 
\\    If $|X_i\cap T| \geq 2$, Add $[B\cup\{v'\},T\setminus \{v\}]$.
\\    If $|X_i\cap T| = 3$, Add $[B\cup \{v''\},T\setminus \{v,v'\}]$.
\\    Similarly for $|X_i\cap T|=4$ and $|X_i\cap T|=5$.\-
\\  If $X_i$ is a  $\overline{K_5}$, For $I=[B,T] \in S$ with $|X_i\cap B| < 2$:\+
\\    If $|X_i\cap B| = 0$:\+
\\      Remove $I$ from $S$; Label elements of $X_i\setminus B$ according to the sequence $v, v', v'', \dots$.
\\      If $|X_i\cap T| \geq 2$, Add $[B\cup\{v,v'\},T]$ to $S$.
\\      If $|X_i\cap T| \geq 3$, Add $[B\cup\{v,v''\},T\setminus \{v'\}]$ and $[B\cup\{v',v''\},T\setminus \{v\}]$
\\      If $|X_i\cap T| \geq 4$, Add $[B\cup\{v,v'''\},T\setminus \{v',v''\}]$, $[B\cup\{v',v'''\},T\setminus \{v,v''\}]$,\+\+\+\+
\\                               and $[B\cup\{v'',v'''\},T\setminus \{v,v'\}]$\-\-\-\-
\\      If $|X_i\cap T| = 5$, Add $[B\cup\{v,v''''\},T\setminus \{v',v'',v'''\}]$, $[B\cup\{v',v''''\},T\setminus \{v,v'',v'''\}]$,\+\+\+\+
\\                                     $[B\cup\{v'',v''''\},T\setminus \{v,v',v'''\}]$, and $[B\cup\{v''',v''''\},T\setminus \{v,v',v''\}]$\-\-\-\-\-
\\    If $|X_i\cap B| = 1$:\+
\\      Similar to if $X_i$ was a $\overline{J_5}$.\-\-\-
\end{algo}
    
For the full details, see our code at \url{https://geometrynyc.wixsite.com/ramsey}.
Once the above process is over, we are left with a set of intervals of potential neighbor sets for $w$.
We enumerate the resulting extended graphs and repeat the above algorithm for each graph. 
Eventually, the process will end for all branches.
We then look for the largest graph that we obtained. 
For the results of this algorithm, see Section \ref{sec:main}.     
    	 
\section{An Alternative Algorithm via a SAT Solver} \label{sec:SAT}
    	 
In this section, we describe an algorithm that handles both gluing and vertex extension.
That is, this algorithm is an alternative to the approach presented in Sections \ref{sec:glue} and \ref{sec:VertexExt}.
One goal of this algorithm is to double check our computations. 
In addition, in some cases this algorithm is faster, partly because it handles the gluing and vertex extensions simultaneously. 

Once again, we follow the notation of Figure \ref{fi:BigPicture}.
In this approach, we turn the problem into a boolean expression and then run a computer program that checks if this expression has a solution. 
We used the CaDiCaL incremental SAT solver.\footnote{\url{https://github.com/arminbiere/cadical}}
For the gluing portion, we create a boolean variable for the existence of every potential edge between $A$ and $B$. 
That is, the edge exists if and only if the variable is true. 
Similarly, for an added vertex $w$, we have a boolean variable for every possible edge between $w$ and another vertex (except $a$ and $b$, which are not connected to additional vertices by definition).

We consider the case where the graph should not contain a $J_6$ and an induced $\overline{K_4}$.
The case of no $K_4$ and no $J_6$ in the complement is handled symmetrically. 
We go over each set of four vertices with at least one from $A$ and at least one from $B$.
We add an or-clause for such a quadruple if every pair of vertices not from $A\times B$ is not connected by an edge.
This clause is false if and only if there are no edges between the four vertices.
In other words, when all boolean variables for the corresponding pairs from $A\times B$ are False. 
These clauses assure that the boolean expression is satisfied only when there is no induced $\overline{K_4}$.

We go over each set of six vertices with at least one from $A$ and at least one from $B$. 
We check how many edges are missing between pairs of vertices not from $A\times B$. 
If exactly one edge is missing, then we create an or-clause that is false if and only if all relevant pairs from $A\times B$ are True.
If zero edges are missing, then we create clauses that are false if and only if at most one edge is missing between these relevant pairs. 
These clauses assure that the boolean expression is satisfied only when there are no copies of $K_6$ and $J_6$.

We first create a boolean formula for the case of a one vertex extension.
If this formula is solvable, then we create a boolean formula for two vertex extensions, and so on. 
The goal is for the process to end before reaching 30 vertices. 
For each new vertex, we add a boolean variable for each potential edge between this vertex and every other vertex.
Similarly to the above, we add clauses for ensuring no $K_6$, $J_6$, and induced $\overline{K_3}$.

To speed up the process we order the new vertices, as follows. 
We represent the set of edges of a vertex as a binary vector and ask these vectors to be ordered (these vectors do not include edges between pairs of added vertices).
Adding such a restriction to the boolean formula by hand is quite difficult. 
Instead, we follow the approach of \cite[Section 3.4]{Zhao17} and use Sympy\footnote{See \url{https://www.sympy.org/en/index.html}.} to convert our expression into the many required clauses.

For the results of this algorithm, see Section \ref{sec:main}.     

\parag{Acknowledgements.} We are grateful for the mentors and organizers of the Polymath Jr program, especially Adam Sheffer, Sherry Sarkar, and David Narvaez.
We thank others from our Polymath Jr Ramsey group for useful conversations, including Mujin Choi, Oliver Kurilov, Nathan Moskowitz, Minh-Quan Vo, Michael Waite, Norbert Weijenberg, and Devin Williams.
Finally, we thank John Mackey for useful discussions.

\section*{Statements \& Declarations.}

\parag{Funding.} This work was supported by NSF award DMS-2218374. 

\parag{Competing interests.} The authors have no relevant financial or non-financial interests to disclose.

\end{document}